\documentclass{amsart} 
\usepackage{amsthm, amsfonts, amssymb, amsmath} 
\newtheorem{theorem}{Theorem}[section]

\theoremstyle{definition}
\newtheorem{definition}[theorem]{Definition}
\newtheorem{example}[theorem]{Example}
\newtheorem{remark/example}[theorem]{Remark/Example}

\newtheorem{question}[theorem]{Question}

 \let\oldlabel=\label
\def\prellabel{\marginparsep=1em\marginparwidth=44pt
 \def\label##1{\oldlabel{##1}\ifmmode\else\ifinner\else
 \marginpar{{\footnotesize\ \\ \tt
  ##1}}\fi\fi}}

\def\cal{}

\def\QQ{ {\bf Q} }

\newcommand{\GL}{\operatorname{GL}}

\newcommand{\Tor}{\operatorname{Tor}}

\newcommand{\PV}{\operatorname{PV}}

\begin{document}

\centerline{ 29th Symposium on Commutative Algebra }
\centerline{ Nagoya November 2007 }

\medskip 

\centerline {\bf \large Koszul algebras and Gr\"obner bases of quadrics}
 
\medskip 

\centerline {Aldo Conca}
\centerline{ Dipartimento di Matematica, Universit\'a di Genova}
\centerline{ Via Dodecaneso 35, I-16146 Genova, Italia}
 
\medskip 

\noindent{ \sl Abstract:}  We present results that appear in the papers  \cite{C,CTV, CRV}  joint with M.E.Rossi, N.V.Trung and G.Valla and also some new results contained in \cite{C1}. These results concern Koszul and G-quadratic properties of algebras associated with points, curves, cubics and spaces of quadrics of low codimension. 
 
\section{Introduction} 

Let $R$ be a standard graded $K$-algebra, that is, an algebra of the form
$R=K[x_1,\dots,x_n]/I$ where $K[x_1,\dots,x_n]$ is a polynomial ring over the
field $K$ and $I$ is a homogeneous ideal with respect to the grading
$\deg(x_i)=1$. Let $M$ be a finitely generated graded $R$-module. 
Consider the (essentially unique) minimal graded $R$-free resolution of $M$ $$\dots\to R^{\beta_i} \to R^{\beta_{i-1}} \to \dots \to R^{\beta_1} \to
R^{\beta_0} \to M\to 0$$ 
The rank $\beta_i$ of the $i$-th
module in the minimal free resolution of $M$ is called the $i$-th {\sl Betti
number} of
$M$. One can also keep track of the graded structure of the resolution. It follows that the free modules in the resolution are indeed
direct sums of ``shifted" copies of $R$: 
$$R^{\beta_i}=\oplus_j R(-j)^{\beta_{ij}}$$ where $R(-a)$ denotes the free
module with the generator in degree $a$ and that the matrices representing the
maps between free modules have homogeneous entries. The number $\beta_{ij}$
is called the {\sl $(i,j)$-th graded Betti number} of $M$. 
 The resolution is finite if $\beta_i=0$ for  $i\gg 0$. 
 
For which algebras $R$ does every module $M$ has a finite minimal free
resolution? The answer is given by (the graded version of) the 
Auslander-Buchsbaum-Serre theorem: 

\begin{theorem}\label{ABS} (Auslander-Buchsbaum-Serre) Let $R$ be a standard graded
$K$-algebra. The following are equivalent: 
 \begin{itemize}
\item[(1)] Every finitely generated graded $R$-module $M$ has a finite minimal free
resolution as an $R$-module. 
\item[(2)] The field $K$, regarded as an $R$-module
via the identification $K=R/\oplus_{i>0} R_i$, has a finite minimal free
resolution as an $R$-module. \par
\item[(3)] $R$ is regular, i.e. $R$ is (isomorphic
to) a polynomial ring.
 \end{itemize}
\end{theorem} 

If $R$ is not regular then the resolution of $K$ is infinite. The {\sl
Poincar\'e series} $P_R(z)$ of $R$ is the formal power series whose
coefficients are the Betti numbers of
$K$, i.e. 
$$P_R(z)=\sum_{i\geq 0} \beta^R_i(K) z^i$$ 
where $\beta^R_i(K)$ is the $i$-th Betti number of $K$ as an $R$-module. 

Serre asked in \cite{S} whether the Poincar\'e series $P_R(z)$ is 
a rational series, that is whether there exist polynomials $a(z),b(z)\in
{\bf Q}[z]$ such that $P_R(z)=a(z)/b(z)$. The positive answer to Serre's
question became well-known under the name of Serre's conjecture. 
This conjecture has been proved for several classes of algebras. For
instance it holds for complete intersections (Tate, Assmus) and for algebras
defined by monomials (Backelin). But in 1981 Anick \cite{A} discovered  
algebras   with irrational Poincar\'e series, such as: 

$$\QQ[x_1,x_2,\dots,x_5]/(x_1^2, x_2^2, x_4^2, x_5^2, x_1x_2, x_4x_5, x_1x_3+x_3x_4+x_2x_5)+m^3$$

 More recently Roos and Sturmfels
have shown that irrational Poincar\'e series arise also in the realm of toric
rings, see\cite{RS}.
 
The Poincar\'e series of $R$ takes into account the rank of the free modules in
the minimal free resolution of $K$. One can also consider the degrees of the
generators of the free modules. This leads to the 
introduction of estimates for the growth of the degrees of the syzygies (like
for instance Backelin's rate) and to the definition of Koszul algebras:

\begin{definition} (Priddy) A standard graded $K$-algebra $R$ 
is Koszul if for all $i$ the generators of the $i$-th free module in the minimal
free resolution of $K$ have degree $i$. Equivalently, $R$ is Koszul if the
entries of matrices representing the maps in the minimal free resolution of
$K$ are homogeneous of degree $1$. 
\end{definition}

\begin{example} Let $R=K[x]/(x^n)$ with $n>1$. Then the resolution of $K$ as an
$R$-module is 
$$ \dots \mathop{\to}^{x^{n-1}} R \mathop{\to}^x R 
 \mathop{\to}^{x^{n-1}} R \mathop{\to}^x R\to K\to 0$$ 
 Hence $R$ is Koszul iff $n=2$. 
 \end{example} 

The algebra $R$ is said to be:
\begin{itemize}
\item[(1)]  quadratic if its defining ideal $I$ is generated by quadrics (i.e. homogeneous elements of degree $2$).  
\item[(2)] G-quadratic if $I$ has a Gr\"obner basis of quadrics with respect to some
system of coordinates and some term order. 
\item[(3)]  LG-quadratic (the L stands for lifting) if there exist a G-quadratic algebra $S$ and a $S$-regular sequence $y_1,\dots, y_s$ of elements of degree $1$ in $S$ such that $S/(y_1,\dots,y_s)\simeq R$. 
\end{itemize} 

One has: 

$$(2)\Rightarrow (3) \Rightarrow \mbox{ Koszul }\Rightarrow (1)$$
 
Implications  $(2)\Rightarrow (3)$  and  $ \mbox{ Koszul }\Rightarrow (1)$ cannot be reversed in general.   We do not know examples of Koszul algebras which are not LG-quadratic. 

By a theorem of Tate (see \cite{F}) every quadratic complete intersection is Koszul, but not all of them are G-quadratic. Non-G-quadratic complete intersection of quadrics are given in \cite{ERT}. The easiest example of non-G-quadratic and quadratic complete intersection is given by $3$ general quadrics in $3$ variables. But every complete intersection of quadrics is LG-quadratic as the following argument of G.Caviglia shows. 

\begin{example} If $R=K[x_1,\dots,x_n]/(q_1,\dots, q_m)$ is a complete intersection of quadrics then $R=S/(y_1,\dots, y_m)$ where $$S=K[x_1,\dots,x_n,y_1,\dots,y_m]/(y_1^2+q_1,\dots, y_m^2+q_m).$$ That $y_1^2+q_1,\dots, y_m^2+q_m$ is a Gr\"obner basis is an easy consequence of Buchberger criterion. That $y_1,\dots, y_m$ form a $S$-regular sequence follows by an Hilbert function computation. 
\end{example}

The Koszul property can be characterized in terms of the Poincar\'e series. Denote
by $H_R(z)$ the Hilbert series of $R$. Then one has:
$$R \hbox{ is Koszul } \Leftrightarrow P_R(z)H_R(-z)=1$$ 
In particular, Koszul algebras have rational Poincar\'e series. 

\section{Filtrations, points and curves}
Given an algebra $R$ it can be very difficult to detect whether $R$ is
Koszul or not. One can  compute the first few matrices in the resolution and check whether they are linear. If they are not, then $R$ is not Koszul. If instead they are linear,  one can then compute few more matrices. But the growth of the size of the matrices (i.e. the growth of the Betti numbers) is in general very fast. And it is known that the first non-linear syzygy can appear in
arbitrarily high homological degree even for algebras with a given Hilbert function.  

A quite efficient method to prove that an algebra is Koszul is that given by filtration arguments of various kinds. These notions have  been used by various authors. Inspired by the work of Eisenbud, Reeves and Totaro \cite{ERT},   Bruns, Herzog and Vetter \cite{BHV} and of Herzog, Hibi and Restuccia \cite{HHR}, we have defined:

\begin{definition} Let $R$ be a standard graded algebra and let $F$ be a 
family of ideals of $R$. Then $F$ is said to be a Koszul filtration of $R$ if
the following conditions hold:
\begin{itemize}
\item[(1)] Every ideal in $F$ is generated by linear
forms,
\item[(2)] The ideal $(0)$ and the maximal homogeneous ideal
$\oplus_{i>0}R_i$ are in $F$,
\item[(3)]For every non-zero $I$ in $F$ there exists
$J$ in $F$ such that
$J\subset I$, $I/J$ is cyclic and $J:I$ is also in $F$. 
\end{itemize} 
\end{definition}

\begin{definition} Let $R$ be a standard graded algebra. A Gr\"obner flag
of $R$ is a Koszul filtration $F$ of $R$ which consists of a single complete
flag. In other words, a Gr\"obner flag is a set of ideals $F=\{(0), (V_1), (V_2),
\dots, (V_n)=(R_1)\}$ where $V_i$ is a $i$-dimensional subspace of $R_1$,
$V_i\subset V_{i+1}$ and $ (V_i):(V_{i+1})=(V_{j})$ for some $j$ depending on $i$. 
\end{definition} 

As the names suggest, we have: 

\begin{theorem} \begin{itemize}
\item[(1)] Let $F$ be a Koszul filtration of $R$. Then
$\Tor^R_i(R/I,K)_j=0$ for all
$i\neq j$ and for all $I\in F$. In particular, $R$ is Koszul. 
\item[(2)] If $R$ has a Gr\"obner flag then $R$ is G-quadratic.
\end{itemize} 
\end{theorem} 

\begin{example}
Let $R=K[x_1,\dots, x_n]/I$ with $I$ a quadratic monomial ideal. 
Then the set $F$ of the ideals of $R$ generated by subsets of $\{x_1,\dots, x_n\}$  is a Koszul filtration of $R$. 
To check it, one has only to observe that  $I:(x_i)$ is generated by variables mod $I$. 
\end{example} 

The property of having a Koszul filtration is stronger than
just being Koszul as the following example shows: 

\begin{example} Let $R$ be a complete intersection of $5$ generic quadrics in $5$. 
As said already above, $R$ is Koszul. But it does not have a Koszul filtration since its defining ideal does not contain quadrics of rank $<3$. 
\end{example} 

As well, to have a Gr\"obner flag is more than G-quadratic. For instance the algebra
$R=K[x,y,z]/(x^2,y^2,xz,yz)$ is obviously G-quadratic but one can easily see
that $R$ does not have a Gr\"obner flag. 

However many classes of algebras which are known to be 
Koszul have indeed a Koszul filtration or even a Gr\"obner flag. For instance:

\begin{theorem} (Kempf) Let $X$ be a set of $s$ (distinct) points of the
projective space ${\bf P}^n$ and let $R(X)$ denote the coordinate ring of
$X$. If $s\leq 2n$ and the points are in general linear position then the
ring $R(X)$ is Koszul, see \cite{K}.
\end{theorem}

We have shown that:

\begin{theorem} With the assumption of Kempf's theorem, 
 the ring $R(X)$ has a Gr\"obner flag.
\end{theorem}

One may ask whether Kempf's theorem holds also for a larger number of points.
This is not the case. 

\begin{example} There exists a set of $9$ points in ${\bf
P}^4$ which are in general linear position and whose coordinate ring is
quadratic but non-Koszul. It is obtained via Gr\"obner-lifting form the ideal number (55) in Roos' list \cite{R}.
\end{example}

 On the other hand for ``generic points" (indeterminates coordinates) we have the following:

\begin{theorem} Let $X$ be a set of ``generic points" in ${\bf P}^n$.
Then $R(X)$ is Koszul if and only if $|X|\leq 1+n+(n^2/4)$.
\end{theorem} 

Let $\cal C$ be a smooth algebraic curve of
genus $g$ over an
algebraically closed field of characteristic zero. If $\cal C$ is not
hyperelliptic, then the
canonical sheaf on $\cal C$ gives a canonical embedding ${\cal C }\to
{\bf P}^{g-1}$ and the
coordinate ring $R_{\cal C}$ of ${\cal C}$ in this embedding is the
canonical ring of $\cal C$.
It is known that $R_{\cal C}$ is quadratic unless $\cal C$ is a trigonal
curve of genus $g\ge 5$ or a plane quintic. 
Another important application of the filtration arguments is the following
theorem. 

\begin{theorem}\label{VF}
 Let $R_C$ be the
coordinate ring of a curve $C$ in its canonical embedding. Assume that $R_C$ is quadratic. 
Then $R_C$ is Koszul. 
\end{theorem}

This is due to Vishik and Finkelberg \cite{VF}; other proofs are given by Polishchuk \cite{P}, and by Pareschi and Purnaprajna \cite{PP}. We are able to show that:

\begin{theorem} Let $R_C$ be as in the Theorem \ref{VF}. Then
$R_C$ has a Gr\"obner flag. 
\end{theorem}

For integers $n,d,s$ the ``pinched Veronese'' $\PV(n,d,s)$ is the $K$-algebra generated by the monomials of degree $d$ in $n$ variables and with at most $s$ non-zero exponents,  that is,    
$$\PV(n,d,s)=K[ x_{1}^{a_1}\cdots x_{n}^{a_n} :\ \ \    \sum a_j=d \ \ \mbox{ and }\ \ \  \#\{ j : a_j>0\}\leq s].
$$   It is an open question whether $\PV(n,d,s)$ is Koszul when quadratic (they are not all quadratic: $\PV(4, 5, 2)$ is not).  G.Caviglia shown in \cite{Ca}  that the first not trivial pinched Veronese $\PV(3,3,2)$ is Koszul by using a combination of filtrations and ad hoc arguments.

\section{Artinian Gorenstein algebras of cubics} 

The algebras $R_C$ are $2$-dimensional Gorenstein domains with h-vector $1+nz+nz^2+z^3$ and Theorem \ref{VF} asserts that they are  Koszul as soon as they are quadratic. 
One might ask: 

\begin{question}
\label{qG1nn1}
Let $R$ be a quadratic Gorenstein algebra  with h-vector $1+nz+nz^2+z^3$.  
Is $R$ Koszul?
\end{question}
Without loss of generality,  one can assume that the algebra  is Artinian. Artinian Gorenstein algebras are described via Macaulay inverse system. 
Let us recall how. Let $S=K[x_1,\dots,x_n]$ be a polynomial ring over a field $K$ of characteristic
$0$. Let $f$ be a non-zero polynomial of $S$ which is homogeneous of degree,
say, $s$. Let $I_f$ be the ideal of $S$ of  the polynomials
$g(x_1,\dots,x_n)$ such that 
$$g(\partial/\partial x_1,\dots,\partial/\partial x_n)f=0. $$ 
Set $R_f=S/I_f$. It is known that $R_f$ is a Gorenstein Artinian
algebra with socle in degree $s$ and that every such an algebra arises in this
way. In particular, in  the case $s=3$  the Hilbert series
of $R_f$ is equal to $1+nz+nz^2+z^3$ (provided $f$ is not a cone).  So Question \ref{qG1nn1} is about algebras   $R_f$ with $f$ a cubic form. 
We are able to show the following:

\begin{theorem} 
\label{balla}
\begin{itemize} 
\item[(1)] Let $f$ be a cubic in $S$. Assume there exist linear forms $y,z$ such that $\partial f/\partial yz=0$ and $\partial f/\partial y$ and $\partial f/\partial z$ are quadrics of rank $n-1$. Then $R_f$ has a Koszul filtration. 
\item[(2)] If $f$ is smooth then $R_f$ is not G-quadratic. 
\item[(3)] For the generic cubic $f$, the ring $R_f$ is Koszul and not G-quadratic. 
\end{itemize} 
\end{theorem} 

Furthermore: 

\begin{theorem} 
\begin{itemize} 
\item[(1)] Let $f$ be a cubic in $S$. Assume there exists linear form $y$ such that $\partial f/\partial y^2=0$ and $\partial f/\partial y$ is quadric of rank $n-1$. Then $R_f$ has a Gr\"obner flag. 
\item[(2)] For the generic singular cubic $f$, the ring $R_f$ is G-quadratic. 
\end{itemize} 
\end{theorem} 

  We are not able to answer Question \ref{qG1nn1} in general.  But in  \cite{CRV} we have shown that Question \ref{qG1nn1} has an affermative answer  $n=3,4$.  In both cases the characterization of the $f$ such that $R_f$ quadratic (or Koszul) is very elegant: 

\begin{theorem}\label{34} For $n=3$ or $4$, the following are equivalent: 
\begin{itemize} 
\item[(1)]  $R_f$ is quadratic. 
\item[(2)]  $R_f$ is Koszul. 
\item[(3)] The ideal of $2$-minors of the Hessian matrix $(\partial f/\partial x_ix_j)$ of $f$ has codimension $n$.
\end{itemize} Furthermore for $n=3$ these conditions are equivalent also to: 
\begin{itemize}
\item[(4)] $f$ is not in the closure of the $\GL_3(K)$-orbit of the Fermat cubic $x_1^3+x_2^3+x_3^3$. 
\end{itemize}
\end{theorem} 

G.Caviglia shown in his unpublished master thesis that property (1),(2) and (3) of \ref{34} are equivalent also in the case of $n=5$.    

Another interesting question is whether the assumption  of \ref{balla}(1) indeed characterize Koszul property for $R_f$. In this case the answer is no, as the following example shows. 

\begin{example} Let $f$ be the Veronese cubic, that is the determinant of a $3\times 3$ symmetric matrix filled with $6$ distinct variables $x_1,\dots,x_6$ and let $H$ be its Hessian matrix. The cubic $f$ has a remarkable property:  $\det H$ is $f^2$ up to scalar and the ideal of $5$-minors of $H$ is $(x_1,\dots,x_6)^2f$. These facts imply that $f$ does not satisfy the assumption of  \ref{balla}(1), nevertheless $R_f$ is Koszul (even G-quadratic).  
\end{example} 
Also, one could also ask whether $R_f$ is LG-quadratic provided it is quadratic. We have reasons to believe that the answer to this question might be positive.

\section{Space of quadrics of low codimension} 
Another point of view we have taken is the following. Let $V$ be a vector space of quadrics of dimension $d$ in $n$  variables. Set $c=\binom{n+1}{2}-d$ the codimension of $V$ in the space of quadrics. Let $R_V$ be the quadratic algebra defined by the ideal generated by $V$. 
A theorem of Backelin \cite{B} asserts that if $c\leq 2$ then $R_V$ is Koszul. We have proved in \cite{C} that: 

\begin{theorem}\label{ciccio}
\begin{itemize}
\item[(1)]  If $c<n$ then  the ring $R_V$ is G-quadratic for  a generic $V$. 
\item[(2)] If $c\leq 2$ then $R_V$ is G-quadratic with, essentially, one exception given by 
$V=\langle  x^2,xy,y^2-xz,yz\rangle$ in $K[x,y,z]$. 
\end{itemize} 
\end{theorem} 

The ``exceptional" algebra  $K[x,y,z]/(x^2,xy,y^2+xz,yz)$ has Hilbert series 
 
$$1+3z+2z^2+z^3+z^4+z^5+\dots$$

and is LG-quadratic since we may deform it to 
$$K[x,y,z,t]/(x^2 + xt,  xy + yt,  yz + xt,  y^2 + xz )$$
which is G-quadratic in the given coordinate system and with respect to revlex  $t>x>y>z$. 

It follows from \ref{ciccio} that every  quadratic Artinian algebra $R$ with  $\dim R_2\leq 2$   is $G$-quadratic. 

A recent conjecture of Polishchuk \cite{P2} on Koszul configurations of points, 
 suggests that   Artinian  quadratic algebras $R$ with  $\dim R_2=3$  should be Koszul. 
This is what we have proved in \cite{C1}: 
 
 \begin{theorem} 
 Let $R$ be an  Artinian  quadratic algebras with  $\dim R_2=3$. Then: 
 \begin{itemize}
 \item[(1)] $R$ is Koszul. 
  \item[(2)]  $R$ is G-quadratic unless it is (up to trivial extension) a complete intersection of $3$ general quadrics in $3$ variables. 
  \end{itemize} 
  \end{theorem}

\end{document}